\documentclass[twoside]{racsam}
\usepackage{amssymb}
\usepackage{graphicx}
\usepackage{amsmath}
\usepackage{amsfonts}
\usepackage{times}

\def\falta{{\bf Falta}}                                         
\def\typeproof{}                                                

\def\commu{\small Art\'{\i}culo panor\'amico / Survey}

\def\area{\small Geometr\'{\i}a y Topolog\'{\i}a / Geometry and Topology}

\def\vol{\bf \falta}
\def\num{\rm \falta}
\def\year{\rm \falta}
\def\presented{\falta.}
\def\submitted{\falta.}
\def\accepted{\falta.}

\begin{document}

\title{Compactness in Banach space theory -- selected problems }

\def\classification{46B26m54D30}
\def\keywords{Hilbertian ball, euclidean ball, Eberlein compact, uniform Eberlein compact, discontinuous norms, weakly compactly generated space, Valdivia compact, Radon-Nikodym compact}


\author{Antonio Avil\'es and Ond\v{r}ej F. K. Kalenda}
\maketitle

\begin{abstract} We list a number of problems in several topics related to compactness in nonseparable Banach spaces. Namely, about the Hilbertian ball in its weak topology, spaces of continuous functions on Eberlein compacta, WCG Banach spaces,  Valdivia compacta and Radon-Nikod\'{y}m compacta. 
\end{abstract}

\begin{titulo}
Compacidad en espacios de Banach - problemas escogidos
\end{titulo}

\begin{resumen} Enumeramos una serie de problemas en diferentes temas relacionados con compacidad en espacios de Banach no separables. Concretamente, sobra la bola eucl\'{\i}dea en su topolog\'{\i}a d\'{e}bil, espacios de funciones continuas en compactos de Eberlein, espacios de Banach d\'{e}bilmente compactamente generados, compactos de Valdivia y compactos de Radon-Nikod\'{y}m.
\end{resumen}


\newtheorem{theorem}{\sffamily Theorem}
\newtheorem{proposition}{\sffamily Proposition}
\newtheorem{lemma}{\sffamily Lemma}
\newtheorem{corollary}{\sffamily Corollary}
\newtheorem{example}{\sffamily Example}
\newtheorem{remark}{\sffamily Remark}
\newtheorem{definition}{\sffamily Definition}
\newtheorem{conjecture}{\sffamily Conjecture}
\def\proof{{\sc Proof.}\quad }
\newtheorem{problem}{\sffamily Problem}
\newtheorem{prop}[theorem]{\sffamily Proposition}
\newtheorem{thm}[theorem]{\sffamily Theorem}
\newtheorem{cor}[theorem]{\sffamily Corollary}

\def\qed{\quad{$\blacksquare$}\medskip }
\def\qedo{\quad{$\square$}\medskip }

\label{firstpage}

\newcommand{\norm}[1]{\left\Vert#1\right\Vert}
\newcommand{\abs}[1]{\left\vert#1\right\vert}
\newcommand{\set}[1]{\left\{#1\right\}}
\newcommand{\Real}{\mathbb R}
\let\R\Real
\newcommand{\eps}{\varepsilon}
\newcommand{\To}{\longrightarrow}
\newcommand{\BX}{\mathbf{B}(X)}
\newcommand{\A}{\mathcal{A}}
\newcommand{\woof}{\mathcal{V}}
\newcommand{\RN}{Radon-Nikod\'ym }
\newcommand{\U}{\mathcal{U}}
\newcommand{\nat}{\mathbb{N}}
\newcommand{\sign}{\mbox{\rm sign\,}}

\section{The hilbertian ball and their relatives}\label{sectionball}

Consider $\kappa$ an uncountable cardinal, identified with the set of ordinals less than $\kappa$. The closed ball of the Hilbert space $\ell_2(\kappa)$ equipped with the weak topology can be identified with the following compact subspace of $\mathbb{R}^\kappa$ (by squaring each coordinate and keeping the sign):

$$B(\kappa) := \left\{x\in\mathbb{R}^\kappa : \sum_{i\in\kappa}|x_i|\leq 1\right\}$$

The closed subspaces of $B(\kappa)$ are the \emph{uniform Eberlein compacta}. This class is well known and studied. Some standard facts are the following:  

\begin{itemize}
\item $B(\kappa)$ is a Fr\'echet-Urysohn space, that is every point in the closure of a set $A\subset B(\kappa)$ is the limit of a sequence of elements of $A$ (this is a consequence of the fact that every element of $B(\kappa)$ has countable support)

\item The density character and the weight of any closed subspace of $B(\kappa)$ coincide. In particular, any separable subspace is metrizable (also a consequence of the fact that every element of $B(\kappa)$ has countable support).

\item Every closed subset of $B(\kappa)$ contains a $G_\delta$ point (take a point with maximal norm).

\item $B(\kappa)$ is sequentially compact: every sequence contains a convergent subsequence (a $G_\delta$-point in the set of cluster points of the sequence provides a convergent subsequence).

\item Every continuous image of a closed subset of $B(\kappa)$ is homeomorphic to a subspace of $B(\kappa)$ \cite{BenRudWag}.
\end{itemize}

These are properties of the class of subspaces of $B(\kappa)$ rather than properties of $B(\kappa)$ itself (actually, these properties are shared by the more general class of Eberlein compacta, though the proofs require much deeper facts). A study of specific topological properties of the space $B(\kappa)$ that distinguish it from its subspaces has been made in recent papers by the authors \cite{unitball,AviKal,KalendaBall,finiteproducts}. The results include ($\kappa$ is an uncountable cardinal):

\begin{itemize}
\item The space $B(\kappa)$ is a continuous image of $A(\kappa)^\mathbb{N}$, where $A(\kappa)$ is the one point compactification of the discrete space of cardinality $\kappa$ \cite{unitball}. It was shown by Bell \cite{Bell} that not all subspaces of $B(\kappa)$ have this property.

\item In \cite{AviKal} a topological invariant called \emph{fiber order} is introduced that --roughly speaking-- associates to a non-metrizable compact set $K$ and a point $x$ a preordered metrizable compact space $\mathbb{F}_x(K)$ and an ordered set $\mathbb{O}_x(K)$ (which is the canonical quotient of $\mathbb{F}_x(K)$). We proved that $\mathbb{O}_x(B(\kappa))$ is a singleton for $x$ in the sphere and $\mathbb{O}_x(B(\kappa))$ is order-isomorphic to the interval $[0,1]$ for other $x$. This allows to check that $B(\kappa)$ is not homeomorphic to its finite powers and to several other spaces.

\item $B(\kappa)$ does not map continuously onto any product $K\times L$ where both $K$ and $L$ are nonmetrizable \cite{finiteproducts}.
\end{itemize}

There is still much to discover about this space. One open problem is the following:

\begin{problem}\label{productwithinterval}
Is $B(\kappa)$ homeomorphic to $B(\kappa)\times K$ for some metrizable compact $K$? For $K = [0,1]$?
\end{problem}

A type of question that can lead the research for understanding this compact space is to know whether it is homeomorphic to some other compact spaces closely related to it. We provide now a family of such spaces. Given two functions $a,b:\kappa\To [-1,1]$ satisfying $-1\leq a(i) < b(i) \leq 1$ for all $i\in\kappa$, we denote
$$B(\kappa,a,b) = \left\{x\in \mathbb{R}^\kappa : \sum_{i\in \kappa}|x_i|\leq 1\text{ and } a(i)\leq x_i \leq b(i) \text{ for all }i\in \kappa\right\}$$

This is a more general family of compact spaces than that consisting of spaces $B(M,N)$ defined in \cite{KalendaBall}, that would correspond to $B(M,a,b)$ where when $b$ is the constant function equal to 1 and $a$ is the function taking the value 0 on $N$ and -1 on $M\setminus N$. Geometrically, for a fixed $i_0\in\kappa$, the set $\{x\in B(\kappa) : a\leq x(i_0)\leq b\}$ corresponds to the intersection of the ball with a \emph{band}, i.e. the region between two hyperplanes, so $B(\kappa,a,b)$ is the intersection of $B(\kappa)$ with a certain family of bands.\\ 

To avoid trivialities, we will assume that 
$$\sum_{i\in\kappa} (a(i)^+ + b(i)^-)<1,$$
where we denote $t^+ = max(0,t)$ and $t^{-} = \max(0,-t)$. We remark that $B(\kappa,a,b)$ is empty if this sum is strictly greater than $1$ and that it is a singleton if the sum equals $1$.

As we shall see later, it is not difficult to check that every space $B(\kappa,a,b)$ is a retract of each other. All the properties of $B(\kappa)$ enumerated above hold for all compacta $B(\kappa,a,b)$, even the computation of the fiber orders $\mathbb{O}_x(K)$ (just the same proof given in \cite{AviKal} is valid for any of these spaces). 

\begin{problem}
Are all compact spaces $B(\kappa,a,b)$ homeomorphic to $B(\kappa)$? Determine the classification up to homeomorphism of the spaces $B(\kappa,a,b)$.
\end{problem}

First we remark that we can restrict ourselves to functions $a$ and $b$ satisfying $-1\leq a(i)\leq 0 <b(i) \leq 1$ and $|a(i)|\leq |b(i)|$ for all $i<\kappa$. We call such pairs \emph{admissible}.

\begin{prop}
Every compact space $B(\kappa,a,b)$ is homeomorphic to a space $(B,\kappa,a',b')$ where $a'$, $b'$ are admissible.
\end{prop}

\proof Suppose that for some $i$ we have either $a(i)>0$ or $b(i)<0$. Denote 
$$A=\{i\in\kappa: a(i)>0\},\qquad B=\{i\in\kappa:b(i)<0\}$$
and set $s=\sum_{i\in A} a(i)-\sum_{i\in B} b(i)$. Then for each $x\in B(\kappa,a,b)$ we have
$$|x_i|\le \begin{cases} 1-s,& i\in\kappa\setminus(A\cup B),\\
1-s+a(i), & i\in A, \\ 1-s-b(i), & i\in B.\end{cases}$$
Therefore we can without loss of generality suppose that
\begin{gather}
|a(i)|\le 1-s,|b(i)|\le 1-s \mbox{ for }i\in\kappa\setminus(A\cup B),\\
b(i)\le 1-s+a(i) \mbox{ for } i\in A,\\
a(i)\ge -(1-s-b(i)) \mbox{ for }i\in B.
\end{gather}
Then the mapping which assigns to each $(x_i)_{i\in\kappa}$ the point $(x'_i)_{i\in \kappa}$ such that
$$x'_i=\begin{cases} \frac{x_i-a(i)}{1-s}, & i\in A, \\ \frac{x_i-b(i)}{1-s}, & i\in B, \\ \frac{x_i}{1-s}, & x\in\kappa\setminus(A\cup B),\end{cases}$$
is a homeomorphism of $B(\kappa,a,b)$ onto $B(\kappa,a',b')$, where
$$a'(i)=\begin{cases} 0 & i\in A, \\ \frac{a(i)-b(i)}{1-s} & i\in B, \\
\frac{a(i)}{1-s}& i\in\kappa\setminus(A\cup B), \end{cases} \qquad 
b'(i)=\begin{cases} \frac{b(i)-a(i)}{1-s} & i\in A, \\ 0 & i\in B, \\
\frac{a(i)}{1-s}& i\in\kappa\setminus(A\cup B). \end{cases}$$
Then $-1\le a'(i)\le 0\le b'(i)\le 1$ for each $i\in\kappa$. Further, for each $i$ such that $|b'(i)|<|a'(i)|$ we can redefine the pair $(a'(i),b'(i))$ by
$(-b'(i),-a'(i))$. In this way we get an admissible pair $(a'',b'')$ such that $B(\kappa,a,b)$ is homeomorphic to $B(\kappa,a'',b'')$.
\qed

 For $t\in [-1,1]$, by abuse of notation, we will denote by $t$ the constant function equal to $t$ and by $t^\ast$ the function on $\kappa$ that it is equal to $t$ except in one point, where it takes value $0$. We can focus on several particular cases:

\begin{itemize}
\item $B(\kappa) = B(\kappa,-1,1)$ is the ball
\item $B^+(\kappa) = B(\kappa,0,1)$ is the positive cone of the ball
\item $B^\frac{1}{2}(\kappa) = B(\kappa,-1^\ast,1)$ is half a ball
\item $B(\kappa, -t,t)$ for $t\in (0,1)$
\item $B(\kappa,0,t)$ for $t\in (0,1)$
\end{itemize}

\begin{problem}\label{fewlistproblem}
Are any two from the compact spaces above homeomorphic?
\end{problem}

This problem is connected with Problem~\ref{productwithinterval}. It is proven in \cite{KalendaBall} that $B^+(\kappa) = B^+(\kappa)\times [0,1]$ and $B^{\frac{1}{2}}(\kappa) = B^{\frac{1}{2}}(\kappa)\times [0,1]$. These results are obtained by using the notions of \emph{cone} $\Delta(K)$ and \emph{double cone} $\Diamond(K)$ of a compact space $K$. The space $\Delta(K)$ is obtained from $K\times [0,1]$ by gluing $K\times \{1\}$ to a singleton, and $\Diamond(K)$ by gluing $K\times\{1\}$ to a singleton and $K\times\{0\}$ to another singleton. Fix $i\in\kappa$ and suppose that the function $a$ satisfies $a(j)\in\{-1,0\}$ for all $j$ and the function $b$ satisfies $b(j)=1$ for all $j\neq i$. We have that:

\begin{itemize}
\item $B(\kappa,a,b) = B(\kappa\setminus\{i\},a,1)\times [0,1]$ if $b(i)<1$
\item $B(\kappa,a,b) = \Delta(B(\kappa\setminus\{i\},a,1))$ if $b(i) = 1$ and $a(i)>-1$;
\item $B(\kappa,a,b) = \Diamond(B(\kappa\setminus\{i\},a,1)$ if $a(i)=-1$.
\end{itemize}

On the other hand, $\Delta(\Delta(K)) = \Diamond(\Delta(K)) =\Delta(\Diamond(K)) = \Delta(K)\times [0,1]$ for any compact space $K$ \cite{KalendaBall}. These facts allow to prove a number a homeomorphic relations between the spaces $B(\kappa,a,b)$ of the form $B(M,N)$ stated in \cite[Theorem 9]{KalendaBall} (the symbol $\sim$ means homeomorphism):
\begin{enumerate}
\item If both $N$ and $M\setminus N$ are infinite, then $B(M,N)\sim \Delta(B(M,N)) \sim \Diamond(B(M,N))\sim B(M,N)\times [0,1]$.
\item If $M\setminus N$ is finite, then $B(M,N)\sim B^+(M) = \Delta(B^+(M)) \sim \Diamond(B^+(M)) \sim B^+(M)\times [0,1]$.
\item If $N$ is finite, then $B(M,N) \sim B^{\frac{1}{2}}(M) \sim \Delta(B^{\frac{1}{2}}(M)) \sim \Diamond(B^{\frac{1}{2}}(M)) \sim B^{\frac{1}{2}}(M) \times [0,1]$.

\end{enumerate}

\begin{prop}\label{retractofeachother}
Each of the spaces $B(\kappa,a,b)$ is a retract of each other.
\end{prop}

\proof It is easy to see that (provided $(a,b)$ is an admissible pair) every $B(\kappa,a,b)$ is a retract of $B(\kappa)$, the retraction sending each $(x_i)_{i<\kappa}$ to $(\min\{\max\{a,x_i\},b\})_{i<\kappa}$. 

Second, we shall prove that $B(\kappa)$ is a retract of $B^+(\kappa)$.
We start by noticing that $B^+(\kappa)$ is homeomorphic to $B^+(\kappa\times\{0,1\})$. We define a continuous surjection of $B^+(\kappa\times\{0,1\})$ onto $B(\kappa)$ by
$$(x_{i,q})_{i<\kappa,q\in\{0,1\}}\mapsto (x_{i,0}-x_{i,1})_{i<\kappa}.$$ 
This surjection has a continuous right inverse $u:B(\kappa)\To B^+(\kappa\times\{0,1\})$  defined by
$$u(x)_{i,0} = x_{i}^+,\qquad u(x)_{i,1} = x_{i}^-,\qquad i<\kappa.$$ 

In order to complete the proof, it remains to show that $B^+(\kappa)$ is a retract of each $B(\kappa,a,b)$. So take $B(\kappa,a,b)$ and consider a partition $\kappa = \bigcup_n J_n$ into infinite or empty sets such that $b(i)>\frac{1}{n}$ if $i\in J_n$. Let $\kappa = \bigcup_{i<\kappa}F_i$ be a partition of $\kappa$ indicated in $\kappa$ with the property that for every $i$ there exists $n$ with $F_i\subset J_n$ and $|F_i|=n$. We define a continuous surjection $f:B(\kappa,a,b)\To B^+(\kappa)$ by 
$$f(x)_i = \sum_{j\in F_i} x_i^+, \qquad i<\kappa.$$
We shall define a  continuous right inverse $u:B^+(\kappa)\To B(\kappa,a,b)$ for $f$. Fix $x\in B^+(\kappa)$ and $i<\kappa$. Enumerate $F_i = \{j_1,\ldots,j_n\}$. We define $u(x)$ on $F_i$ in the following way:
$$u(x)_{j_1} = \min\{x_i,b(j_1)\},
u(x)_{j_2} = \min\{x_i-u(x)_{j_1},b(j_2)\},\ldots,u(x)_{j_n} = \min 
\{x_i - \sum_{l=1}^{n-1}u (x)_{j_l},b(j_n)\}.$$ 
\qed

We shall focus now in one particular instance of Problem~\ref{fewlistproblem}:

\begin{problem}
Are $B(\kappa)$ and $B^+(\kappa)$ homeomorphic?
\end{problem}

The nice thing about this \emph{nonseparable} problem is that is equivalent to a \emph{separable} problem. In the following theorem $\|x\| := \sum_{i\in\kappa} |x_i|$ (notice that, if we identify $B(\kappa)$ with the ball of $\ell_2(\kappa)$ by means of the homomorphism described above, this is the square
of the Euclidean norm; for this theorem that square is irrelevant). Also, for $x\in [-1,1]^{\aleph_0}$, $\|x\|_\infty = \sup\{|x_n| : n\in\aleph_0\}$.

\begin{thm}\label{equivalentproblems} The following statements are equivalent.
\begin{enumerate}
\item For every uncountable cardinal $\kappa$, $B(\kappa)$ is homeomorphic to $B^+(\kappa)$.
\item There is an uncountable cardinal $\kappa$ such that $B(\kappa)$ is homeomorphic to $B^+(\kappa)$.
\item There exists a homeomorphism $f:B(\aleph_0)\To B^+(\aleph_0)$ such that $\|f(x)\| < \|f(y)\|$ whenever $x,y\in B(\aleph_0)$, $\|x\| < \|y\|$.
\item There exists a homeomorphism $f:B(\aleph_0)\To B^+(\aleph_0)$ such that $\|f(x)\| = \|x\|$ for all $x\in B(\aleph_0)$
\item There exists a homeomorphism $f:[-1,1]^{\aleph_0}\To [0,1]^{\aleph_0}$ in the product topologies, such that $\|f(x)\|_\infty = \|x\|_\infty$ for all $x$.
\end{enumerate}
\end{thm}

\proof $(1\Rightarrow 2)$ is evident.

$(2\Rightarrow 3)$ requires using the technique of fiber orders from~\cite{AviKal}. As it is computed there, for $x\in K$ the fiber preorder structure $\mathbb{F}_x(B(\kappa))$ is either trivial (a singleton) or it is the compact space $B(\aleph_0)$ endowed with the preorder relation ($x\leq y \iff \|x\|\leq \|y\|$). Just the same computation shows that $\mathbb{F}_x(B^+(\kappa))$ is either trivial or it is the compact space $B^+(\aleph_0)$ endowed with the order ($x\leq y \iff \|x\|\leq \|y\|$). Hence, if $B(\kappa)$ and $B^+(\kappa)$ are homeomorphic, then there must exist an order-preserving homeomorphism between nontrivial fiber preorders $\mathbb{F}_x(B(\kappa))$ and $\mathbb{F}_y(B^+(\kappa))$.

For $(3\Rightarrow 4)$ it is enough to prove the following statement:

Claim: For every strictly increasing bijection $\phi:[0,1]\To [0,1]$ there exists a homeomorphism $g:B^+(\aleph_0)\To B^+(\aleph_0)$ such that $\|g(x)\| = \phi(\|x\|)$.

The implication follows from the claim because if $(3)$ holds, then there is a strictly increasing bijection $\psi:[0,1]\To [0,1]$ such that $\|f(x)\| = \psi(\|x\|)$; apply the claim to $\phi = \psi^{-1}$ and consider $g\circ f$. For proving the claim we define $g(x_0,x_1,x_2,\ldots) = (y_0,y_1,y_2,\ldots)$ as follows: $y_0 = \phi(x_0)$, $y_n = \phi(\sum_{0}^n x_i) - \sum_0^{n-1}y_i$.

Finally, for $(4\Rightarrow 1)$, consider a partition of $\kappa$ into countable sets. For each of these countable sets $S$ we have a norm-preserving homeomorphism $B(S)\To B^+(S)$. Gluing all these homeomorphisms together we get a (norm-preserving also) homeomorphism $B(\kappa)\To B^+(\kappa)$

It remains to prove the equivalence of $4$ and $5$. This follows from the following theorem, which is a consequence of results of Dijkstra and van Mill \cite{DijMil}.

\begin{thm}\label{dijmil} There exists a homeomorphism $f:[-1,1]^{\aleph_0}\To B(\aleph_0)$ such that $f([0,1]^{\aleph_0})=B^+(\aleph_0)$ and $\|f(x)\| = \|x\|_\infty$ for every $x\in [-1,1]^{\aleph_0}$.
\end{thm}

\proof It follows from \cite[Theorem 5.2]{DijMil} that there is an onto homeomorphism $g:\R^{\aleph_0}\To\R^{\aleph_0}$ such that:

\begin{eqnarray}
\label{DM1} \forall x\in \R^{\aleph_0}: \|g(x)\|_1=\|x\|_{\infty} \qquad \mbox{(including the infinite values)} \\
\label{DM2} g\mbox{ is sign-preserving, i.e. }\forall x\in \R^{\aleph_0}\forall n\in\aleph_0: g(x)_n=0 \mbox{ or }\sign g(x)_n=\sign x_n
\end{eqnarray}

Let $f$ be the restriction of $g$ to $[-1,1]^{\aleph_0}$. Then $f$ is a norm-preserving homeomorphism
of $[-1,1]^{\aleph_0}$ onto $B(\aleph_0)$. This follows immediately from (\ref{DM1}). 
It remains to show that $f([0,1]^{\aleph_0})=B^+(\aleph_0)$.

First notice that it follows from (\ref{DM2}) that for every $x\in\mathbb{R}^{\aleph_0}$ and every $n\in\aleph_0$:

\begin{equation}
		\label{DM3} (x_n\ge 0\Rightarrow g(x)_n\ge0) \text{ and } (x_n=0\Rightarrow g(x)_n=0).
\end{equation}

Hence it is clear that $f([0,1]^{\aleph_0})\subset B^+(\aleph_0)$. It remains to prove the reverse inclusion. As finitely supported vectors are dense in $B^+(\aleph_0)$ and $f([0,1]^{\aleph_0})$ is compact, it is enough to show that the range contains all finitely supported vectors in $B^+(\aleph_0)$.
To prove it we introduce the following notation. If $F\subset\aleph_0$ is finite, we set
\begin{align*}
U(F)&=\{x\in [0,1]^{\aleph_0} : x_i=0\mbox{ for }i\notin F\},\\
V(F)&=\{x\in B^+({\aleph_0}) : x_i=0\mbox{ for }i\notin F\}
\end{align*}
It follows from (\ref{DM1}) and (\ref{DM3}) that $f(U(F))\subset V(F)$ for each finite $F\subset\aleph_0$. If we show that $f(U(F))=V(F)$ for each $F$, the proof will be completed.

We will show it by induction on the cardinality of $F$. If $F$ is a singleton, then it is obvious.
Indeed, if $e_i$ denotes the $i$-th canonical vector in $R^{\aleph_0}$, it follows from (\ref{DM1}) and (\ref{DM3}) that $g(te_i)=te_i$ for each $t\in\R$.

Further, let $n\in\nat$ be such that the equality holds for each $F$ of cardinality at most $n$. Suppose that $F$ has cardinality $n+1$. Fix $t\in(0,1]$ and set
\begin{alignat*}{3}
C&= \{x\in U(F): \|x\|_\infty=t\},&  \qquad C_i&=\{x\in C: x_i=0\}\mbox{ for }i\in F,\\
D&= \{x\in V(F): \|x\|_1=t\},& D_i&=\{x\in D: x_i=0\}\mbox{ for }i\in F.
\end{alignat*}
Then clearly $f(C)\subset D$. Moreover, by the induction hypothesis $f(C_i)=D_i$ for each $i\in F$.
We observe that $C$ is homeomorphic to $B_n$, the $n$-dimensional euclidean ball. Indeed, the function $\phi(x) = (x_i - \frac{1}{n+1}\sum_{j\in F}x_j)_{i\in F}$ provides a homeomorphism from $C$ onto the ball of radius $t$ of the space $\{u\in\mathbb{R}^F : \sum_{j\in F}u_j = 0\}$ endowed with the norm $\|u\|' = \sup\{|u_i-u_j| : i\neq j\}$. On the other hand, the pair $(D,\bigcup_{i\in F} D_i)$ is homeomorphic to the pair $(B_n,S_{n-1})$, where $S_{n-1}$ denotes the sphere of $B_n$. This is because $\Delta = D\setminus\bigcup_{i\in F}D_i$ is a convex open subset of the affine space $\{x\in\mathbb{R}^{\aleph_0} : x_i=0\mbox{ for } i\not\in F,\ \sum_{j\in F}x_j = t\}$ and $D=\overline{\Delta}$. Finally, $f(C)$ is a subset of $D$ which is homeomorphic to $D$ and contains $\bigcup_{i\in F}D_i$. Thus $f(C)=D$ (as the only subset of $B_n$ which is homeomorphic to $B_n$ and contains $S_{n-1}$ is the whole set $B_n$). 

As $t\in(0,1]$ was arbitrary, we get $f(U(F))=V(F)$ which completes the proof. 
\qed

\section{Spaces of continuous functions on some Eberlein compacta}\label{sectionCK}

In the preceeding section we dealt with the topological structure of $B(\kappa)$, and now we shall deal with the structure of the Banach space of continuous functions $C(B(\kappa))$. The first remark is the following:

\begin{prop}\label{CBkisomorphic}
Each of the Banach spaces $C(B(\kappa,a,b))$ is isomorphic to $C(B(\kappa))$.
\end{prop}

In order to prove that, one uses Pe\l czy\' nski's decomposition method:

\begin{prop}[{\rm Pe\l czy\'nski}]
Let $X$ and $Y$ be Banach spaces. Suppose that $X$ is isomorphic to a complemented subspace of $Y$ and $Y$ is isomorphic to a complemented subspace of $X$. Suppose also that $X$ is isomorphic to $c_0(X)$, the $c_0$-sum of contably many copies of $X$. Then $X$ and $Y$ are isomorphic.\\
\end{prop}

Proof of Proposition \ref{CBkisomorphic}: If $K$ is a retract of a compact space $L$, then $C(K)$ is a complemented subspace of $C(L)$, so from Proposition~\ref{retractofeachother} we conclude that the spaces $C(B(\kappa,a,b))$ are complemented in each other. A result of Argyros and Arvanitakis \cite{ArgArv} that $c_0(C(B(\kappa)))$ is isomorphic to $C(B(\kappa))$ completes the proof.\qed

The situation is nevertheless unclear to us for some compact spaces that we know that are not homeomorphic, but we do not know about the spaces of continuous functions:

\begin{problem}\label{CBkn}
Is $C(B(\kappa))$ isomorphic to $C(B(\kappa)^n)$ for $n\leq\omega$? Is $C(B(\kappa))$ isomorphic to $C(P(K))$, where $P(K)$ is any of the spaces of probability measures considered in \cite{AviKal}: $P(A(\kappa)^n)$, $P(\sigma_n(K))$...?
\end{problem}

Motivated by the classical Milutin's theorem, it was an open problem for long time whether every Banach space of continuous functions $C(K)$ is isomorphic to a space $C(L)$ with $L$ a zero-dimensional compactum. A counterexample to this question was first found by Koszmider \cite{Koszmider}. The key property of this example, and those constructed later in a similar spirit, is that it is an \emph{indecomposable Banach space}, i.e., it contains no nontrivial complemented subspaces. However all the Banach spaces we are considering here are weakly compactly generated, hence they are very rich in complemented subspaces. Therefore, the problem is still worth being investigated under restrictions (like being WCG) that make indecomposable counterexamples impossible. 

\begin{problem}\label{CBk}
Is $C(B(\kappa))$ isomorphic to some $C(L)$ with $L$ a zero-dimensional compact space?
\end{problem}

Even the following weaker version of the problem stated by Argyros and Arvanitakis is open \cite{ArgArv}:

\begin{problem} Is there a nonseparable weakly compact convex set $K$ in a Banach space and a zero-dimensional compact space $L$ such that $C(K)\cong C(L)$?
\end{problem}

Hence Problem \ref{CBk} is open also for any of the spaces mentioned in Problem~\ref{CBkn}, and in particular for $K = B(\kappa)^n$ with $n$ finite or countable. Argyros and Arvanitakis \cite{ArgArv} were able to construct a weakly compact star-shaped set $K$ containing a ball such that $C(K)$ is isomorphic to $C(L)$ with $L$ zero-dimensional. Now, we can propose a specific candidate for $L$. If $C(L)$ should be isomorphic to $C(B(\kappa))$, $L$ has to be uniform Eberlein and not scattered. We propose $L = A(\kappa)^\mathbb{N}$, where $A(\kappa)$ is the one-point compactification of a discrete set of cardinality $\kappa$.\\

\begin{problem}
Is $C(B(\kappa))$ or some $C(B(\kappa)^n)$ isomorphic to $C(A(\kappa)^\mathbb{N})$?\\
\end{problem}

Let us have a look at what we know concerning the conditions of Pelczynski's decomposition method when dealing with these problems. First thing is to know if $X\cong c_0(X)$. In the litterature we find:

\begin{itemize}
\item $C(A(\kappa)^\mathbb{N})$ is isomorphic to $c_0(C(A(\kappa)^\mathbb{N})$ \cite{countableproducts}.
\item $C(B(\kappa))$ is isomorphic to $c_0(C(B(\kappa))$ \cite{ArgArv}.\\
\end{itemize}

The second thing is to know which spaces are isomorphic to complemented subspaces of other spaces. In the case of spaces of continuous functions $C(K)$, the usual tool is that of averaging and extension operators. An \emph{extension operator} for a continuous injection $v:K\To L$ is an operator $T:C(K)\To C(L)$ such that $T(f)(v(x)) = f(x)$ for all $x\in K$ and $f\in C(K)$. An \emph{averaging operator} for a continuous surjection $u:L\To K$ is an operator $T:C(L)\To C(K)$ such that $T(u\circ f) = f$ for all $f\in C(K)$. In either of the two case, the existence of such an operator $T$ implies that $C(K)$ is complemented in $C(L)$. Often one finds averaging or extension operators that are regular, that is: positive, of norm one, and sending constant functions to constant functions of the same value. Let us see a few remarks about what this technique can say and what it cannot say for our purposes.

\begin{prop}$C(A(\kappa)^\mathbb{N})$ is isomorphic to a complemented subspace of $C(B(\kappa)^\mathbb{N})$.\end{prop}

\proof Consider the the canonical embedding $u:A(\kappa)\To B(\kappa)$ that sends every element $i\in\kappa$ to the $x\in B(\kappa)$ that is zero everywhere except at $i$ where $x_i=1$. This embedding has a regular extension operator: To every continuous function $f\in C(A(\kappa))$ we can associate its extension $T(f)\in C(B(\kappa))$ defined as $T(f)(x) = f(\infty)\left(1-\sum_{i\in\kappa}|x_i|\right) + \sum_{i\in\kappa} |x_i|f(i)$. By \cite[Proposition 4.7]{Pelczynski}, the countable power $u^\mathbb{N}:A(\kappa)^\mathbb{N}\To B(\kappa)^\mathbb{N}$ also admits a regular extension operator, hence $C(A(\kappa)^\mathbb{N})$ is isomorphic to a complemented subspace of $C(B(\kappa)^\mathbb{N})$.\qed

\begin{prop}
If $m<n\leq\omega$, then $C(B(\kappa)^m)$ is a complemented subspace of $C(B(\kappa)^n)$.
\end{prop}

\proof $B(\kappa)^m$ is a retract of $B(\kappa)^n$, and this provides an embedding with regular extension operator. \qed

The problem with the powers of $B(\kappa)$ arises in finding $C(B(\kappa)^n)$ complemented in $C(B(\kappa)^m)$ for $m<n$. We cannot hope to get this using averaging operators since $B(\kappa)^m$ does not map continuously onto $B(\kappa)^n$ by \cite{finiteproducts}. Next proposition shows that neither regular extension operators can be used.

\begin{prop}
If $m<n\leq\omega$, there does not exist any embedding $B(\kappa)^n\To B(\kappa)^m$ with a regular extension operator.
\end{prop}

\proof Suppose $u:B(\kappa)^n\To B(\kappa)^m$ is such an embedding. A regular extension operator would provide a continuous function $g: B(\kappa)^m\To P(B(\kappa)^n)$ with $gu(x) = \delta_x$ for $x\in B(\kappa)^n$ \cite{Pelczynski}. Since $B(\kappa)^n$ can be viewed as a compact convex set, we can consider the barycenter map $h:P(B(\kappa)^n)\To B(\kappa)^n$ which is continuous, affine, and satisfies $h(\delta_x) = x$. The composition $hg:B(\kappa)^m\To B(\kappa)^n$ is continuous and onto, since $hg(u(x))= x$. However such onto maps do not exist for $m<n$ by \cite{finiteproducts}.\qed

\begin{prop}
The continuous surjection $f:A(\kappa)^\mathbb{N}\To B(\kappa)$ constructed in \cite{unitball} does not admit any averaging operator.
\end{prop}

Proof: Let $r_n$ a sequence of positive real numbers that defines a continuous surjection $2^\mathbb{N}\To [0,1]$ by $(t_n)\mapsto \sum r_n t_n$. For example $r_n = 2^{-n-1}$. Let $Z = \{z\in\mathbb{N}^\mathbb{N} : \sum_{n=0}^\infty \frac{z_n}{r_n} \leq 1\}$ and $L = Z\times \prod_{n\in\mathbb{N}}\sigma_n(\kappa)^\mathbb{N}$. Denote elements of $L$ as $(z,x)$, where $x = (x[m,n]\in\sigma_n(\kappa) : n,m\in\mathbb{N}$). A continuous surjection $g:L \To B^+(\kappa)$ is defined as $g(z,x)_i = \sum_{m\in\mathbb{N}}\frac{1}{r_m} x[m,z_m]_i$. The surjection $f$ from \cite{unitball} is the composition of $g$ with a continuous surjection $g:A(\kappa)^\mathbb{N}\To L$. It is enough to see that $g$ has no averaging operator. This can be seen using a result of Ditor. Define a derivation procedure on $B^+(\kappa)$ in the following way: given a set $X\subset B^+(\kappa)$, its derived set $X'$ is the set of all $x\in X$ whose fiber $g^{-1}(x)$ contains two disjoint sets $S$ and $T$ that belong to the closure in the Vietoris topology of the set of fibers $\{g^{-1}(y) : y\in X\}$. Inductively, $X^{(n)} = (X^{(n-1)})'$. According to \cite[Corollary 5.4]{Ditor}, if $B^+(\kappa)^{(n)} \neq\emptyset$ for all $n$, then $g$ does not admit any averaging operator. Let $A = \{y\in B^{+}(\kappa) : \sum_{i\in\kappa}|y_i|<1\}$. We prove that $A' = A$. Let $y\in A$, and $\varepsilon = 1-\sum_{i\in\kappa}|y_i|$. We can consider a sequence of elements $y_n\in A$ converging to $x$ in which $y_n$ coincides with $y$ on the support of $y$, and $y_n$ has a single extra nonzero coordinate with value $\frac{3}{4}\varepsilon$. It is possible to find a natural number $k$ big enough so that for every $(z,x)\in g^{-1}(y_n)$ we have $\sum_{j<k}z_j \geq 1-\varepsilon/2$. Now let $m$ another integer big enough so that $\frac{3}{4m}\varepsilon<\sum_{j\geq k}r_j$. We pick a new sequence $y'_n\in A$ converging to $y$, in which $y'_n$ coincides with $y$ on the support of $y$ and has exactly $m$ extra nonzero coordinates with value $\frac{3}{4m}\varepsilon$ on each of these coordinates. The elements $y'_n$ have the property that for every $(z,x)\in g^{-1}(y'_n)$, $\sum_{i<k}z_i\leq 1-\frac{3}{4}$. Two cluster points in the Vietoris topology $S$ and $T$ of the sequences $\{g^{-1}(y_n)\}$ and $\{g^{-1}(y'_n)\}$ respectively are two disjoint closed sets of $g^{-1}(y)$, since $S\subset\{(z,x) : \sum_{j<k}z_j \geq 1 -\varepsilon/2\}$, while $T\subset\{(z,x) : \sum_{j<k}z_j \leq 1-\frac{3}{4}\varepsilon\}$. This proves that $y\in A'$.\qed

\begin{problem}
Is there a continuous surjection $f:A(\kappa)^\mathbb{N}\To B(\kappa)$ with an averaging operator?
\end{problem}

\section{Eberlein compacta, WCG spaces and their subspaces}

There are certain natural hierarchies of classes of Banach spaces. Two of them --the so called \emph{descriptive hierarchy} and \emph{differentiability hierarchy}-- are described for example in \cite[Chapter 1]{habil}. Many of these classes can be characterized using a topological property of the weak topology. For example we have:

\begin{thm} Let $X$ be a Banach space.
\begin{itemize}
	\item $X$ is separable if and only if $(X,w)$ is separable.
	\item $X$ is reflexive if and only if $(X,w)$ is $\sigma$-compact.
	\item $X$ is WCG if and only if $(X,w)$ contains a dense $\sigma$-compact subset.
	\item $X$ is weakly Lindel\"of determined (WLD) if and only if $(X,w)$ is primarily Lindel\"of.
	\item $X$ is Asplund if and only if each separable subset of $(X,w)$ can be covered by countably many relatively closed metrizable subsets.
\end{itemize}
\end{thm}
 
There exist a number of characterizations of subspaces of WCG spaces (cf.~\cite{FMZ} and \cite[Theorem 6.13]{biorthogonal}), but all of them involve the linear structure of the Banach space.

\begin{problem} Is it possible to characterize subspaces of WCG spaces in terms of the weak topology? More precisely, suppose that $X$ is a subspace of a WCG space, and $Y$ is Banach space weakly homeomorphic to $X$, is $Y$ also a subspace of a WCG space?
\end{problem}

This problem has is meaningful as subspaces of WCG spaces need not be WCG. There are two basic examples -- the first one, by Rosenthal \cite{rosenthal}, is a subspace of a large $L^1(\mu)$ space; the second one is due to Argyros \cite[Theorem 1.6.3]{Fabiansbook} and it is a natural subspace of a $C(K)$ space with $K$ uniform Eberlein compact. But there are many spaces which are hereditarily WCG, apart from separable or reflexive spaces, an example is $c_0(\Gamma)$ for an arbitrary set $\Gamma$. So, the following problem seems to be natural. 

\begin{problem} Characterize those compact spaces $K$ for which $C(K)$ is hereditarily WCG.
\end{problem}

This question was posed to us by Marian Fabian some time ago. Scattered Eberlein compacta have this property. Indeed, if $K$ is scattered and Eberlein,
then $C(K)$ is simultaneously Asplund and WCG. Further, WCG Asplund spaces are hereditarily WCG by \cite[Corollary 6]{OriVal}, cf. also \cite[Proposition 8.3.2]{Fabiansbook}. Another class with the property that $C(K)$ is hereditarily WCG is that of Eberlein compacta of weight less than $\mathfrak b$ by \cite[Corollary 7]{RNcompact}. Cardinal $\mathfrak b$ is the \emph{bounding number}: the least cardinality of a subset of $\mathbb{N}^{\mathbb{N}}$ which is not contained in any $\sigma$-compact subset of $\mathbb{N}^{\mathbb{N}}$. Consistently $\mathfrak b >\aleph_1$. On the other hand, the above mentioned example $C(K)$ of Argyros can be taken so that $K$ is a uniform Eberlein compact space of weight $\mathfrak b$. Hence $C(B(\mathfrak b))$ and $C(A(\mathfrak b)^\mathbb{N})$ are not hereditarily WCG, since every uniform Eberlein compact space of weight $\kappa$ is a subspace of $B(\kappa)$ and a continuous image of a subspace of $A(\kappa)^\mathbb{N}$. 

One more problem inspired by the fact that WCG is not a hereditary property is the following one:

\begin{problem} Let $K$ be a convex Eberlein compact set, i.e. a convex compact subset of a locally convex space which is Eberlein compact. Is it homeomorphic
to a convex compact subset of $(X,w)$ for a Banach space $X$?
\end{problem}

Recall that one of the equivalent definitions of an Eberlein compact space is that these are just compact subsets of $(X,w)$ for a Banach space $X$. If $X$ is a subspace of WCG, then $(B_{X^*},w^*)$ is convex Eberlein compact (the underlying locally convex space being $(X^*,w^*)$). But $(B_{X^*},w^*)$ is \emph{affinely} homeomorphic to a convex compact subset of $(Y,w)$ for a Banach space $Y$ if and only if $X$ is WCG.

Indeed, suppose first that $X$ is WCG. Let $K\subset (X,w)$ be a generating compact subset. We can embed $(B_{X^*},w^*)$ affinely homeomorphically onto a bounded subset of $C_p(K)$ using the canonical restriction map $x^*\mapsto x^*|K$. As bounded pointwise compact subsets of $C(K)$ are weakly compact, the first implication is proved.

Conversely, suppose that $(B_{X^*},w^*)$ is affinely homeomorphic to a convex compact set $K\subset (Y,w)$ for a Banach space $Y$. Note that without loss of generality we can suppose that the spaces $X$ and $Y$ are real. (If they are complex, we consider them over the field of reals. This preserves all the respective properties.) Consider the standard restriction map defined by $y^*\mapsto y^*|K$. This is a continuous map of $(B_{Y^*},w^*)$ into $C_p(K)$. Moreover, the image $L$ of $B_{Y^*}$ separates points of $K$. This is a standard way to prove that $C(K)$ is WCG. But in our case, we suppose that $K$ is convex and so, $L\subset A(K)$, the Banach space of all affine continuous functions on $K$. It is now easy to check that $L\cup\{1\}$ is a weak compact generating subset of $A(K)$. Indeed, every element $\phi\in A(K)^*$ is of the form
$$\phi(f)=a f(x) - b f(y),\quad f\in A(K)$$
for some $x,y\in K$ and $a,b\in\mathbb{R}$. (To see it extend $\phi$ to an element of $C(K)^*$. This extension is of the form $a\mu_1-b\mu_2$ for some $a,b\in\mathbb{R}$ and Radon probabilities $\mu_1$,$\mu_2$. Then set $x$ to be the barycenter of $\mu_1$ and $y$ to be the barycenter of $\mu_2$.)
So, $L\cup\{1\}$ is a weak compact subset of $A(K)$ which contains $1$ and separates points of $K$. It follows from the above described representation of the dual that the only $\phi\in A(K)^*$ which is zero on $L\cup\{1\}$ is the zero functional.
So, by Hahn-Banach theorem $L\cup\{1\}$ generates $A(K)$, hence $A(K)$ is WCG.
But $A(K)$ is isometric to $A(B_{X^*},w^*)$ which is isomorpic to $X\times\mathbb{R}$. Thus $X$ is WCG as a hyperplane of a WCG space.

So, there are examples of Banach spaces $X$ such that $(B_{X*},w*)$ is Eberlein but it cannot be affinely embedded into $(Y,w)$ for any Banach space $Y$. But it is not clear whether it must be homeomorphic to a convex subset of some $(Y,w)$.

\section{Valdivia compacta and their interactions with other structures}

Valdivia compacta form a rich and natural class of compact spaces which is useful in the investigation of nonseparable Banach spaces with the focus on decompositions to separable pieces. For a detailed study of this class we refer to \cite{survey}. A collection of natural examples of Valdivia compacta can be found in \cite{valexa}, and a corresponding study on a more general class in \cite{kubisjmaa}.
We recall the necessary definitions.

\let\S\Sigma  \let\G\Gamma
Let $\Gamma$ be an arbitrary set. We define
$$\S(\G)=\{ x\in\R^\G : \{\gamma\in\G : x(\gamma)\ne 0\}
\mbox{ is countable} \}.$$
Further, let $K$ be a compact space.
\begin{itemize}
\item We
say that $A\subset K$ is a {\it $\S$-subset of $K$} if there is a
homeomorphic injection $h$ of $K$ into some $\R^\G$ such that
$A=h^{-1}(\S(\G))$.
\item $K$ is called a {\it Valdivia compact space}
if $K$ has a dense $\S$-subset.
\end{itemize}

Moreover, $K$ is called a {\it Corson compact space}
if $K$ is a $\S$-subset of itself, i.e. if it is homeomorphic to a subset of $\S(\G)$.

We start with the following problem. A subset $A$ of a complex vector space is called absolutely convex if $\sum_{i=1}^n t_i x_i\in A$ whenever $x_1,\ldots,x_n\in A$ and $t_1,\ldots,t_n\in\mathbb{C}$ with $\sum_{i=1}^n|t_i|\leq 1$.

\begin{problem} Let $X$ be a complex Banach space. Suppose that $(B_{X^*},w^*)$ admits a convex symmetric dense $\S$-subset. Does it admit an absolutely convex one?
\end{problem}

The answer is positive for $C(K)$ spaces (see \cite[Theorem 3.8]{complex}) and in case $(B_{X^*},w^*)$ admits only one dense $\S$-subset (which is trivial).
A related problem is the following one.

\begin{problem} Let $K$ be a Valdivia compact space and let $G$ be a compact (or even finite) abelian group of automorphisms of $K$. Is there a dense $\Sigma$-subset of $K$ which is invariant for each element of $G$?
\end{problem}

Even the following simplest case seems to be open.

\begin{problem} Let $K$ be a Valdivia compact space and $h:K\to K$ a homeomorphism satisfying $h\circ h=\mbox{id}$. Is there an $h$-invariant dense $\Sigma$-subset of $K$?
\end{problem}

The following problem is inspired by results of \cite{nova}.

\begin{problem} Is there a Radon-Nikod\'ym Valdivia compact space which contains a copy of $[0,\omega_2]$? Or even a scattered one?
\end{problem}

We recall that a Valdivia compact space which does not contain a copy of $[0,\omega_1]$ is Corson \cite{survey} and that a Corson Radon-Nikod\'ym compact space is Eberlein \cite{OriSchVal}. A positive answer would yield essentially new examples of Valdivia compacta, a negative one would yield that the dual of an Asplund space with countably $1$-norming Markushevich basis also contains such a basis (see \cite[Section 6]{nova}).

Another question from \cite{nova} is the following.

\begin{problem} Let $X$ be an Asplund space. Is there an equivalent norm on $X$ such that the respective bidual unit ball is Valdivia compact (in the weak* topology)?
Or even such that the dual admits a countably $1$-norming Markushevich basis?
\end{problem}

If $X$ is Asplund and has density $\aleph_1$, then $(B_{X^{**}},w^*)$ is Valdivia compact by \cite[Corollary 4.5]{nova}. In this case no renorming is needed.
On the other hand, there is an Asplund space of density $\aleph_2$ such that the respective bidual unit ball is not Valdivia compact\cite[Example 4.10(b)]{nova}.
This space is a renorming of $C[0,\omega_2]$. However, if we consider this space
with the canonical norm, the bidual unit ball is Valdivia compact (as it is for each $C(K)$ space, see e.g. \cite[Example 4.10(a)]{nova}). In this context it is worth to mention the result of \cite[Proposition 26]{kubisjmaa}, which implies that the dual to an Asplund space admits a $1$-projectional skeleton, hence the bidual unit ball admits a retractional skeleton. The definitions of these notions can be found in \cite{kubisjmaa}. We only stress that this is a kind of a noncommutative analogue of Valdivia compacta (Valdivia compact spaces are exactly those admitting a commutative retractional skeleton by \cite{kubmich}).
So, the question is whether we can by a renorming ensure commutativity.

The next problem is from \cite{valexa} and \cite{valdline}.

\begin{problem}\label{orderpreservingimage} Let $K$ be a linearly ordered compact space which is a continuous image of a Valdivia compactum. Is $K$ an order-preserving continuous 
image of a linearly ordered Valdivia compactum?
\end{problem}

Linearly ordered Valdivia compacta have a quite nice structure by \cite[Section 3]{valexa}. A characterization of such spaces is the content of \cite{valdline}.
The structure of linearly ordered compact spaces which are continuous images of Valdivia compacta is less clear. A criterion for recognizing whether $K$ is an order-preserving continuous image of a linearly ordered Valdivia compactum is contained in \cite[Theorem 6.1]{valdline}. In view of this the following question
is also interesting.

\begin{problem}\label{uncountablecharacter} Let $K$ be a linearly ordered compact space which is a continuous image of a Valdivia compactum. Is the cardinality of the set of all points of uncountable character at most $\aleph_1$?
\end{problem}

A positive answer to Problem~\ref{orderpreservingimage} would imply a positive answer to Problem~\ref{uncountablecharacter}, due to \cite[Proposition 5.5]{K2006a} and \cite[Theorem 6.1]{valdline}.
However, the opposite direction is not clear due to \cite[Example 3.5]{valdline}.

The last problem is from \cite[Section 7]{valexa}.

\begin{problem} Let $A$ be a $C^*$ algebra. Is the bidual unit ball $(B_{A^{**}},w^*)$ Valdivia compact? Is it true for $A=B(\ell_2)$, the algebra of all bounded linear operators on $\ell_2$?
\end{problem}

For commutative $C^*$ algebras, i.e. for spaces $C_0(T)$ with $T$ locally compact, the answer is positive (see e.g. \cite[Theorem 5.5]{valexa}).
More generally, the answer is positive for type I $C^*$ algebras, i.e. for those $C^*$ algebras for which the bidual is a von Neumann algebra of type I (which is in this case equivalent to semifiniteness). This follows from \cite[Theorem 7.1]{valexa}.

\section{Radon-Nikod\'{y}m compact spaces}

Radon-Nikod\'{y}m (RN) compact spaces are defined as weak$^\ast$ compact subsets of duals of Asplund spaces. Ever since it was asked by Namioka \cite{Namioka}, the main open question on this subject has been whether the continuous image of an RN compactum is RN compact. 

\begin{problem}
Is the continuous image of an RN compact space also RN compact?
\end{problem}

There has been a number of papers addressing this problem \cite{OriSchVal,StegallRN,Arvanitakis, FabHeiMat,MatSte,Namiokanote,RNcompact,RNorder,ArvAvi}, including a survey article \cite{Fabiansurvey}. We collect in this section a list of problems scattered in the litterature and some remarks, all of them related to the continuous image problem. First, we recall the intrinsic topological characterization of RN compacta. Let $d:K\times K\To [0,+\infty)$ be a function with the property that $d(x,y)=0$ if and only if $x=y$. The function $d$ is said to \emph{fragment} $K$ if for every nonempty $L\subset K$ and every $\varepsilon>0$ there exists $U\subset K$ open such that $U\cap L\neq\emptyset$ and $\sup\{d(x,y) : x,y\in U\cap L\}<\varepsilon$. A compact space is $RN$ if and only if there is a lower semicontinuous metric $d$ on $K$ that fragments $K$. If one does not assume $d$ to satisfy the triangle inequality one obtains the concept of \emph{quasi $RN$ compact}. The class of quasi RN compacta is closed under continuous images \cite{Arvanitakis} but the following is unknown:

\begin{problem}
Is every quasi RN compact space an RN compactum? Is every quasi RN compact space a continuous image of an RN compactum?
\end{problem}

The class of quasi RN compacta was discovered independently by different authors \cite{Arvanitakis, FabHeiMat}, and it admits several nice characterizations \cite{RNcompact, Fabiansurvey}. If $K$ is quasi $RN$ compact, then $K$ is $RN$ compact if any of the following two conditions hold: (1) the weight of $K$ is less than $\mathfrak b$ \cite{RNcompact}, or (2) $K$ is \emph{almost totally disconnected} \cite{Arvanitakis}, i.e., $K\subset [0,1]^I$ and for every $x\in K$, $|\{i\in I : x_i\not\in\{0,1\}\}|\leq\aleph_0$ (that includes $K$ being Corson, Eberlein or zero-dimensional).

All classes of quasi RN, RN and their continuous images are closed under the following operations: passing to subspaces, making countable products, taking spaces of probability measures $P(K)$, and taking spaces of Borel measures of variation at most $1$, $M_1(K) = (B_{C(K)^\ast},w^\ast)$. The closed convex hull of compact space $K$ is a continuous image of the space $P(K)$, so the following problem by Namioka is a particular case of the continuous image problem:

\begin{problem}
Let $K$ be a compact subset of a locally convex space. Assume that $K$ is RN compact. Is the closed convex hull of $K$ also RN compact?
\end{problem}

Arvanitakis \cite{Arvanitakis} proved the following:

\begin{thm}\label{0dimRN} For a zero-dimensional compact space, the following are equivalent
\begin{enumerate}
\item $K$ is RN
\item $K$ is quasi RN
\item $K$ embeds in a countable product of scattered compact spaces.
\end{enumerate}
\end{thm}

In this direction, Argyros asked whether the whole class of RN compacta (and not only zero-dimensional) can be generated in a sense from the class of scattered compacta. Namely, 

\begin{problem}
Let $K$ be an RN compact space. Is $K$ homeomorphic to a subspace of $P(S)$ for $S$ scattered compact? Is it a subspace of $M_1(S)$? Is it a continuous image of a subspace of $M_1(S)$?
\end{problem}

\begin{problem}
Let $K$ be a RN compact (or quasi RN compact) space. Is $K$ the continuous image of a zero-dimensional RN compact space $L$? Can a continuous surjection $f:L\To K$ be obtained with an averaging operator?  
\end{problem}

The interest of getting the averaging operator is the relation with the continuous image problem. If $f:L\To L$ has an averaging operator, then $C(L)$ is a complemented subspace of $C(K)$. A compact space is RN compact if and only if $C(K)$ is an Asplund generated space \cite[Theorem 1.5.4]{Fabiansbook}, and this class is closed under complemented subspaces. The fact that these classes of compact spaces can be characterized through their $C(K)$ Banach spaces, allows to establish a connection with the problem of getting $C(K)$ isomorphic to $C(L)$ with $L$ zero-dimensional, mentioned in Section~\ref{sectionCK}.

\begin{prop}
Suppose that there exists a compact space $K$ that is quasi RN compact but not RN compact. Then $C(K)$ is not isomorphic to $C(L)$ for any zero-dimensional compact space $L$.
\end{prop}

Proof: Suppose that $C(K)$ is isomorphic $C(L)$. Since $K$ is quasi RN compact space, $C(K)$ is a $\sigma$-Asplund generated Banach space \cite{sigmaasplund}. Hence so is $C(L)$, so $L$ is quasi RN compact \cite{sigmaasplund}. But $L$ being zero-dimensional, Arvanitakis' Theorem \ref{0dimRN} implies that $L$ is RN. Hence $C(L)$ is isomorphic to $C(K)$ is Asplund generated, and  $K$ is RN \cite[Theorem 1.5.5]{Fabiansbook}. \qed

The approach from \cite{FabHeiMat} to quasi RN compact spaces and the recent characterization of continuous images of RN compacta in \cite{ArvAvi} suggest a relevant role of the combinatorics of $\mathbb{N}^\mathbb{N}$ in the problem of distinguishing quasi RN, RN compacta and their continuous images. We have also the fact that the three classes coincide for compact spaces of weight less than $\mathfrak b$ \cite{RNcompact}. Nevertheless, we do not know the answer to the following

\begin{problem}
Assume that every continuous image of an RN compactum $K$ of weight $\mathfrak c$ is RN compact. Does this imply that every continuous image of an RN compact space is RN compact? Suppose that all quasi RN compacta of weight $\mathfrak c$ are continuous images of RN compacta. Does this imply that this holds for all quasi RN compacta?
\end{problem}

A natural way of addressing this problem is through this other one:

\begin{problem}
Let $K$ be a  continuous image of an RN compact space. Suppose that all continuous images of $K$ of weight at most $\mathfrak c$ are RN compact. Is $K$ RN compact? Suppose that $K$ is qRN and all its continuous images of weight $\mathfrak c$ are continuous images of RN compacta. Is $K$ a continuous image of an RN compactum? 
\end{problem}

Of course, we may ask by the way a similar question like:

\begin{problem}
Let $K$ be a  continuous image of an RN compact space. Suppose that all closed subspaces of $K$ of weight $\leq \mathfrak c$ are RN compact. Is $K$ RN compact? Similarly for quasi RN compacta and continuous images of RN compacta.
\end{problem}

After the results from \cite{RNcompact}, it would be natural to ask the same problems for cardinal $\mathfrak b$ instead of $\mathfrak c$, because quasi RN compact spaces of weight less than $\mathfrak b$ are RN, and it is for weight $\mathfrak b$ where first essential difficulties in the continuous image problem appear.

In \cite{ArvAvi}, Arvanitakis and the first author developped a method of construction of continuous images of RN compacta that are natural candidates not to be RN compact. We mention here a concrete one: the space of almost increasing functions. Let $(\prec)$ be a well order on $\mathbb{N}^\mathbb{N}$, and $$AIF := \left\{x\in [0,1]^{\mathbb{N}^\mathbb{N}} : x_\sigma \leq x_\tau + \frac{1}{2^{min\{n : \sigma_n\neq\tau_n\}}}\text{ whenever }\sigma\prec\tau\right\}.$$
The space $AIF$ is a continuous image of an RN compactum \cite{ArvAvi}.

\begin{problem}
Is $AIF$ an RN compact space?
\end{problem}

Finally we mention another particular instance of the problem of continuous images studied by Matou\v{s}kov\'{a} and Stegall \cite{MatSte}:

\begin{problem}
Let $K$ be a compact space that can be written as the union of two RN compact spaces $K_1$ and $K_2$. Is $K$ RN compact?
\end{problem}

Notice that the union of two compact spaces is the continuous image of the disjoint union, and clearly the dijoint union of two RN compact spaces is RN compact. Matou\v{s}kov\'{a} and Stegall \cite{MatSte} prove that the answer to this question is positive in some particular cases: when $K_1\cap K_2$ is either metrizable, scattered or a $G_\delta$ set, or when $K_1$ is a retract of $K$ or when $K\setminus K_1$ is scattered. Fabian \cite{Fabiansurvey} asks whether one could prove the result for $K_1\cap K_2$ Eberlein compact.

%



\begin{acknow}
The first author was supported by MEC and FEDER (Project MTM2008-05396), Fundaci\'{o}n S\'{e}neca (Project 08848/PI/08) and Ramon y Cajal contract (RYC-2008-02051).

The second author was supported in part by the grant GAAV IAA 100190901 and in part by the
Research Project MSM 0021620839 from the Czech Ministry of Education.
\end{acknow}

Antonio Avil\'{e}s, Departamento de Matem\'{a}ticas, Universidad de Murcia, Campus de Espinardo, 30100 Murcia (Spain).  

Ond\v{r}ej F.K. Kalenda, Department of Mathematical Analysis, 
Faculty of Mathematics and Physic, Charles
University, Sokolovsk\'a 83, 186 75 Praha 8, Czech Republic.
\label{lastpage}
\end{document}